\documentclass[dvips,nixaihp,preprint]{imsart}
\RequirePackage[OT1]{fontenc}
\RequirePackage{amsthm,amsmath,natbib,amssymb}
\RequirePackage[colorlinks,citecolor=blue,urlcolor=blue]{hyperref}

\startlocaldefs
 \numberwithin{equation}{section}
 \theoremstyle{plain}
 \newtheorem{thm}{Theorem}[section]
 
 \newtheorem{cor}[thm]{Corollary}
 
 \newtheorem{lem}[thm]{Lemma}
 
 \newtheorem{lemma}[thm]{Lemma}
 \newtheorem{prop}[thm]{Proposition}
 \theoremstyle{definition}

 \theoremstyle{remark}
 
 \newtheorem{rem}[thm]{Remark}
 
 \newtheorem{example}{Example}
 \numberwithin{equation}{section}
 \newcommand{\corr}{\operatorname*{cor}}

 \newcommand{\R}{\mathbb{R}}

 \newcommand{\Z}{\mathbb{Z}}

 \newcommand{\E}{\mathbb{E}}
 \renewcommand{\P}{\mathbb{P}}

\endlocaldefs

\begin{document}

\begin{frontmatter}
\title{
Intermittency and Aging for the Symbiotic Branching Model}
\runtitle{Intermittency and Aging for Symbiotic Branching}

\begin{aug}
\author{\fnms{Frank} \snm{Aurzada}\ead[label=e1]{aurzada@math.tu-berlin.de}},
\author{\fnms{Leif} \snm{D\" oring}\thanksref{t2}\ead[label=e2]{doering@math.tu-berlin.de}}
\thankstext{t2}{Supported by the DFG International Research Training Group ``Stochastic models of Complex Processes'', Berlin.}

\runauthor{Aurzada, D\" oring}

\affiliation{Technische Universit\"at Berlin\thanksmark{m1}}

\address{Institut f\"ur Mathematik\\
Technische Universit\"at Berlin\\
Stra\ss e des 17.~Juni 136\\ 10623 Berlin\\ Germany\\
~\\
\phantom{E-mail:\ }\printead*{e1}
\phantom{E-mail:\ }\printead*{e2}}

\end{aug}

\begin{abstract}
	For the symbiotic branching model introduced in \cite{EF04}, it is shown that aging and inter\-mit\-tency exhibit different behaviour for negative, zero, and positive correlations. Our approach also provides an alternative, elementary proof and refinements of classical results concerning second moments of the para\-bolic Anderson model with Brownian potential. Some refinements to more general (also infinite range) kernels of recent aging results of \cite{DD07} for interacting diffusions are given.
\end{abstract}

\begin{keyword}[class=AMS]
\kwd[Primary ]{60K35}
\kwd[; secondary ]{60J55}
\end{keyword}

\begin{keyword}
\kwd{Aging}
\kwd{Interacting Diffusions}
\kwd{Intermittency}
\kwd{Mutually Catalytic Branching Model}
\kwd{Parabolic Anderson Model}
\kwd{Symbiotic Branching Model}
\end{keyword}

\end{frontmatter}


\section{Introduction}
For the last three decades, equations of the type 
\begin{eqnarray}
	d u(t,i)=\sum_{j\in\Z^d}a(i,j)(u(t,j)-u(t,i))\,dt+\sqrt{\kappa f(u(t,i))}\,dW_t(i)\label{111}
\end{eqnarray}
have been studied intensively. Here, $i\in \Z^d$, $t\geq 0$, $\kappa>0$, $(a(i,j))_{i,j\in\Z^d}$ are transition rates on $\Z^d$, and $W=\{W_t(i)|t\geq 0, i\in\Z^d\}$ is a familiy of independent Brownian motions. The following special cases with very different interpretations and different behaviour are quite common in the literature.

\begin{example}\label{ex1}
	The (Wright-Fisher) stepping stone model from mathematical genetics: $f(x)= x(1-x)$.
\end{example}

\begin{example}\label{ex2}
	The parabolic Anderson model (with Brownian potential) from mathematical physics: $f(x)= x^2$.
\end{example}

\begin{example} \label{ex3}
	The super random walk from pure probability theory: $f(x)= x$.
\end{example}

\begin{example}\label{ex4}
	The critical (spatial) Ornstein-Uhlenbeck process: $f(x)=1$.
\end{example}
For the super random walk, $\kappa$ is the branching rate which in this case is time-space independent. In \cite{DP98}, a two type model based on two super random walks with time-space dependent branching was introduced. The branching rate for one species is proportional to the value of the other species. More precisely, the authors considered 
\begin{align*}
	d u(t,i)&=\sum_{j\in\Z^d}a(i,j)(u(t,j)-u(t,i))\,dt+\sqrt{\kappa u(t,i)v(t,i)}\,dW^1_t(i),\\
	d v(t,i)&=\sum_{j\in\Z^d}a(i,j)(v(t,j)-v(t,i))\,dt+\sqrt{\kappa u(t,i)v(t,i)}\,dW^2_t(i),
\end{align*} 
where now $W^1=\{W^1_t(i)|t\geq 0, i\in\Z^d\}$ and $W^2=\{W^2_t(i)|t\geq 0, i\in\Z^d\}$ are families of independent Brownian motions. Solutions are called mutually catalytic branching processes. In the following years, properties of this model were well studied (see for instance \cite{CK00}, \cite{CDG04}).\\

In this paper, we are interested in a variant of mutually catalytic branching, namely the symbiotic branching model introduced in \cite{EF04} for continuous space. The same equations as for the mutually catalytic branching model are considered but additionally the driving noises are correlated in the following way:
\begin{eqnarray}
	\langle W_{\cdot}^n(i),W_{\cdot}^m(j)\rangle_t \, = \, 
       	\begin{cases}
		\varrho t&:i=j\text{ and }n\neq m,\\
		t&: i=j\text{ and } n=m,\\
		0&: \text{otherwise},
	\end{cases}\label{22} 
\end{eqnarray}
where $\varrho\in[-1,1]$ is a correlation parameter. For $\varrho=0$ and as well for general $\varrho$ there are basically two approaches to formalize the equations. In \cite{DP98} under quite restrictive assumptions on the transition kernel $(a(i,j))_{i,j\in\Z^d}$, existence of solutions was obtained in the space of tempered sequences. Since their assumptions in particular assume symmetry and exponential decay of $(a(i,j))_{i,j\in\Z^d}$, already existence of solutions in cases we are interested in is not clear. This is why we stick to the setup of \cite{CDG04}, which as well is more popular for interacting diffusions. For the transition kernel $(a(i,j))_{i,j\in\Z^d}$ we assume
\begin{align*}
	&(H_1)\,\,\,\, 0\leq a(i,j)<\infty,\\
	&(H_2)\,\,\,\, \sum_{j\in\Z^d}a(i,j)=1,\quad \forall i\in\Z^d,\\
	&(H_3)\,\,\,\, a(i,j)=a(0,i-j).	
\end{align*}
Two main examples of interest are the following:
\begin{example}\label{e1}
	The discrete Laplacian is given by
	\begin{eqnarray*}
		a(i,j)=\begin{cases}
		        \frac{1}{2d}&:|i-j|=1,\\
		        0&:\textrm{otherwise}.
		       \end{cases}
	\end{eqnarray*}
	Obviously, $(H_1)$, $(H_2)$, $(H_3)$ are fulfilled. Further, the one-dimensional Riemann walk (see for instance \cite{H95}) has transition rates
	\begin{eqnarray*}
		a(i,j)=a(0,|i-j|)=  \frac{c}{|i-j|^{1+\beta}},
	\end{eqnarray*}
	with $c$ normalising the total rate to $1$. Here $(H_1)$, $(H_2)$, $(H_3)$ are also fulfilled but in contrast to the discrete Laplacian the assumptions of \cite{DP98} are not satisfied.
\end{example}

To specify the state space we fix a positive, summable function $\alpha$ on $\Z^d$ satisfying 
\begin{eqnarray*}
	\sum_{i\in\Z^d}\alpha(i)a(i,j)\leq K\alpha(j),\quad\forall j\in\Z^d,
\end{eqnarray*}
for some finite constant $K$. In \cite{CDG04} a possible choice (see their Equation (4.13)) is given by
\begin{eqnarray*}
	\alpha(j)=\sum_{i\in\Z^d}\sum_{n=0}^{\infty}\frac{1}{K^n}p^{(n)}(i,j)\beta(i),
\end{eqnarray*}
where $\beta$ is positive and summable, $p^{(n)}$ denote the $n$-step transition probabilities, and $K>1$. 
This is needed to verify the generalized Mytnik self-duality which was introduced for the continuous space analogue model in Proposition 5 of \cite{EF04}. The duality for the discrete space is similar. For the duality in $E$ in the special case $\varrho=0$ see Lemma 4.1 of \cite{CDG04}. The state space is now defined by pairs of functions of the following Liggett-Spitzer space
\begin{eqnarray} 
	E=\big\{f:\Z^d\rightarrow \R_{\geq 0}\big|\sum_{j\in\Z^d} f(j)\alpha(j)<\infty\big\}.
\end{eqnarray}
The choice of $\alpha$ does not influence the results.
\begin{prop}\label{o}
	For $u_0,v_0\in E$ and $\varrho\in[-1,1]$ there is a (weak) solution of the symbiotic branching model with almost surely continuous paths and state space $E$. For $\varrho\in[-1,1)$ solutions are unique in law.
\end{prop}
The proof of Proposition \ref{o} is standard. Existence can be proven by finite dimensional approximations as in \cite{SS80}. For $\varrho\in(-1,1)$ uniqueness follows from the generalized Mytnik self-duality as in \cite{EF04}. For $\varrho=-1$ uniqueness is true since moments increase slowly enough, and for $\varrho=1$ uniqueness of  solutions is not known.\\ 
For this work we restrict ourselves to homogeneous initial conditions
\begin{eqnarray*}
	\,\,\,\,u_0=v_0\equiv 1.
\end{eqnarray*}
	This is not necessary but simplifies the notation a lot.\\
	
	The interesting feature of the symbiotic branching model is that it connects Examples~\ref{ex1}-\ref{ex3} above. Being first established as a time-space inhomogeneous version of a pair of super random walks, the examples from above appear as special cases: $\varrho=0$ obviously corresponds to the mutually catalytic branching model. The case $\varrho=-1$ with the additional assumption $u_0+v_0\equiv 1$ corresponds to the stepping stone model as can be seen as follows: Since in the perfectly negatively correlated case $W^1(i)=-W^2(i)$, the sum $u+v$ solves a discrete heat equation and with the further assumption $u_0+v_0\equiv 1$ stays constant for all time. Hence, for all $t\geq 0$, $u(t,\cdot)\equiv 1-v(t,\cdot)$ which shows that $u$ is a solution of the stepping stone model with initial condition $u_0$ and $v$ is a solution with initial condition $v_0$. Finally, suppose $w$ is a solution of the parabolic Anderson model, then, for $\varrho=1$, the pair $(u,v):=(w,w)$ is a solution of the symbiotic branching model with initial conditions $u_0=v_0=w_0$.\\

	The purpuse of this and the accompanying paper \cite{BDE09} is to understand the nature of the symbiotic branching model better. How does the model depend on the correlation $\varrho$? Are properties of the extremal cases $\varrho\in\{-1,0,1\}$ inherited by some regions of the parameters? Since the longtime behaviour of super random walk, stepping stone model, mutually catalytic branching model, and parabolic Anderson model are very different, one might guess that for varying $\varrho$ different regimes correspond to the different models.\\

	The focus of \cite{BDE09} lies on the longtime behaviour in law (unifying the classical results for the stepping stone model, mutually catalytic branching model, and parabolic Anderson model) if $(a(i,j))_{i,j\in\Z^d}$ generates a recurrent Markov process, on (un)boundedness of higher moments $\E[u(t,k)^p]$ as $t\rightarrow \infty$, and the wave speed for the continuous space analogue. It was shown that in the recurrent case $\E[u(t,k)^p]$ is bounded in $t$ if and only if 
	\begin{eqnarray}\label{e}
		p<\frac{\pi}{\pi/2+\arctan\big(\varrho/\sqrt{1-\varrho^2}\big)}.
	\end{eqnarray}
	For the transient case the behaviour in $\varrho$ is open to a large extent.\\
	
	In contrast to \cite{BDE09} the present paper focuses on second moments. Note that (\ref{e}) implies that in the recurrent case second moments are bounded if and only if $\varrho<0$. This can also be seen and analysed in more detail using a moment-duality which will be explained in Section \ref{sm}. Using this duality, we show how to reduce second moments of symbiotic branching processes to moment generating functions and Laplace transforms of local times of discrete space Markov processes, i.e.
	\begin{eqnarray*}
		\E[e^{\kappa L_t}],
	\end{eqnarray*}
	where $\kappa\in\R$ and $L_t$ denotes the local time $\int_0^t\delta_0(X_s)\,ds$ in $0$. For simple random walks the behaviour as $t\to\infty$ was partially analysed in \cite{CM94} by analytic methods. Here, we present a simple new proof based on a renewal-type equation and Tauberian theorems. The simplicity of the proof has the advantage that no further assumptions on the Markov process (in particular no symmetry and no finite range assumptions) are needed. For any $\kappa>0$ and $\kappa<0$ the technique yields precise growth rates including all constants.\\
	As an application, intermittency and aging for symbiotic branching processes are established. \\	

	The main results on intermittency and aging are collected in Section~\ref{sec:results}. In Section \ref{section2} we first establish representations for second moments of symbiotic branching processes (Section \ref{sm}), then prove the results for exponential moments of local times (Section \ref{sectionlocal}), and, finally, proofs of the main results (Section \ref{sectionmain}) are given.
\section{Results} \label{sec:results}

\subsection{Intermittency}
The first property we address is intermittency (see for instance \cite{CM94} for a discussion of the ideas). The $p$-th Lyapunov exponent is defined by
\begin{eqnarray}
	\gamma^u_p(\kappa,\varrho)=\lim_{t\rightarrow \infty}\frac{1}{t}\log \E[u(t,k)^p]\label{ly}
\end{eqnarray}
if the limit exists ($\gamma^v_p(\kappa,\varrho)$ analogously). Since in this work we only deal with second moments we further define
\begin{eqnarray*}
	\gamma_2^{u,v}(\kappa,\varrho)=\lim_{t\rightarrow \infty}\frac{1}{t}\log \E[u(t,k)v(t,k)].
\end{eqnarray*}
In Lemma~\ref{duality} we will see that $\gamma_2^{u,v}(\kappa,\varrho)=\gamma_2^{u}(\kappa,\varrho)=\gamma_2^{v}(\kappa,\varrho)$ and hence we abbreviate $\gamma_2$. 
One says the system is intermittent (or weakly intermittent as recently in \cite{FK09}) if $\gamma_2>0$. 

Intermittency for the parabolic Anderson model ($\varrho=1$) is a well-studied property (see \cite{CM94} and \cite{GdH07}). The existing proofs heavily depend on the linear structure of the system since they employ explicit solutions given as Feynman-Kac type representations. Such explicit solutions are not known to exist for the symbiotic branching model. Hence, one might ask whether or not the results obtained for $\varrho=1$ can be transferred to some larger regime of correlation values. Indeed, this can be done.\\
Let us first fix some notation. In the following $(X_t)$ denotes a continuous time Markov process with transition rates $(a(i,j))_{i,j\in\Z^d}$ and $p_t(i,j)=\P[X_t=j|X_0=i]$.
Due to the moment-duality for symbiotic branching processes (see Lemma \ref{lem:moment-duality}) the notation of symmetrization is needed. For two independent Markov processes  $(X^1_t)$, $(X^2_t)$ with transition rates $(a(i,j))_{i,j\in\Z^d}$ the symmetrization is defined as 
\begin{eqnarray}\label{sy}
	\bar X_t=X^1_t-X^2_t.
\end{eqnarray}
The transition rates of the symmetrization are given by
\begin{eqnarray*}
	\bar a(i,j)=a(i,j)+a(j,i),
\end{eqnarray*}
its transition probabilities are denoted $\bar p_t(i,j)$. Note that in the symmetric case $\bar p_t(i,j)=p_{2t}(i,j)$. The Green function of $(X_t)$ is denoted $G_{\infty}(i,j)=\int_0^{\infty} p_t(i,j)\,dt$ and we abbreviate $G_{\infty}= G_{\infty}(0,0)$. Further, we set $H_{\infty}(i,j)=\int_0^{\infty}t\,p_t(i,j)\,dt$ and abbreviate $H_{\infty}=H_{\infty}(0,0)$. Analogously, we use $\bar G_{\infty}, \bar H_{\infty}$ for the symmetrization $(\bar X_t)$.
\begin{thm}[Weak Intermittency for Symbiotic Branching]\label{t3}Let $(u_t,v_t)$ be a solution of the symbiotic branching model with homogeneous initial conditions. Then $(u_t,v_t)$ is intermittent if and only if
			\begin{eqnarray*}
				\kappa\varrho>\frac{1}{ \bar G_{\infty}}.
			\end{eqnarray*}
	In particular, there is no intermittency for non-positive $\varrho$.
\end{thm}

The previous theorem suggests $0$ deviding symbiotic branching into two regimes in which the $\varrho>0$ regime behaves like the parabolic Anderson model with respect to intermittency.\\
Although after understanding the two-types particle moment-dual (Lemma \ref{mom}) one sees that the problem can be treated as for $\varrho=1$, we present a new proof. In \cite{GdH07} results of \cite{CM94} for higher moments of the parabolic Anderson model were generalized to more general symmetric transitions than the discrete Laplacian. Here, in particular, complete results for the asymptotic behaviour (exponential and subexponential) of second moments of the parabolic Anderson model with arbitrary transitions $(a(i,j))_{i,j\in\Z^d}$ are proven.\\

In the course of the proofs we obtain the expression
\begin{eqnarray}
	\gamma_2=\hat{\bar p}^{-1}\Big(\frac{1}{\kappa\varrho}\Big),
\end{eqnarray}
for the second Lyapunov exponent, where $\hat{\bar p}^{-1}$ is the inverse of the Laplace transform of the return probabilities (see Proposition \ref{prop1}). This expression, by Tauberian theorems, gives us the explicit asymptotic behaviour for the Lyapunov exponents as function of $\kappa$. In the following, $\sim$ denotes strong asymptotic equivalence, i.e.\ $h_1\sim h_2$ means $\lim h_1/h_2=1$.

\begin{prop}\label{cor}Let $(u_t,v_t)$ be a solution of the symbiotic branching model with homogeneous initial conditions. Then for $\kappa\varrho>\frac{1}{\bar G_{\infty}}$, the map $\gamma_2:\big[\frac{1}{\varrho \bar G_{\infty}},\infty\big)\rightarrow \R_{\geq 0} ,\kappa\mapsto\gamma_2(\kappa,\varrho)$ has the following properties:
	\begin{itemize}
 		\item[i)] $\gamma_2$ is strictly convex,
 		\item[ii)] $\gamma_2(\kappa)\leq \kappa\varrho$ for all $\kappa$, and $\frac{\gamma_2(\kappa)}{\kappa\varrho} \to 1$ for $\kappa\to\infty$,
 		\item[iii)] if $\bar p_t(0,0)\sim c t^{-\alpha}$ as $t\to \infty$, $\alpha\leq 1$, we have, as $\kappa\to 0$,
			\begin{align*}
			\gamma_2(\kappa) \sim \begin{cases}
   			(c \Gamma(1-\alpha) \kappa\varrho)^{1/(1-\alpha)} &: 0<\alpha<1,\\
   			\exp( -(c \kappa\varrho)^{-1} + o(\kappa^{-1})) &: \alpha=1,\\
			\end{cases}
			\end{align*}
		\item[iv)] if $\bar p_t(0,0)\sim c t^{-\alpha}$, as $t\to \infty$, $\alpha>1$, we have, as $\kappa\searrow \frac{1}{(\bar G_{\infty}\varrho)}>0$,
			\begin{align*}
			\gamma_2(\kappa) \sim \begin{cases}
   			\left(\frac{(\kappa\varrho-1/\bar G_{\infty}) \bar G_{\infty}^2 (\alpha-1) }{c \Gamma(2-\alpha)} \right)^{1/(\alpha-1)} &: 1<\alpha<2,\\
   			\frac{\bar G_{\infty}^2}{c} (\kappa\varrho-1/\bar G_{\infty}) (\log 1/(\kappa\varrho-1/\bar G_{\infty}))^{-1} &: \alpha=2,\\
   			\frac{\bar G_{\infty}^2}{\bar H_{\infty}} (\kappa\varrho-1/\bar G_{\infty}) &: \alpha>2.
			\end{cases}
			\end{align*}
	\end{itemize}
	Here, $\Gamma$ denotes the Gamma function. 
\end{prop}
Our approach has the further advantage that the growth rates in the critical and subcritical regimes follow directly:
\begin{prop}\label{pro}
	Let $(u_t,v_t)$ be a solution of the symbiotic branching model with homogeneous initial conditions. 
	If $\bar p_t(0,0)\sim c t^{-\alpha}$, as $t\rightarrow \infty$, then the following hold:
	\begin{itemize}
	 \item $\varrho>0$ and $\alpha>1$
	\begin{itemize}
		\item[i)] If $\kappa\varrho<\frac{1}{\bar G_{\infty}}$, then
			\begin{eqnarray*}
				\E[u(t,k)^2]\sim \frac{1}{\varrho(1-\kappa\varrho \bar G_{\infty})},\quad\text{ as }t\to \infty.
			\end{eqnarray*}
		\item[ii)] If $\kappa\varrho=\frac{1}{\bar G_{\infty}}$ then, as $t\rightarrow \infty$,
			\begin{eqnarray*}
			\E[u(t,k)^2] \sim \begin{cases}
   			\frac{\bar G_{\infty}(\alpha-1)}{c \Gamma(2-\alpha)\Gamma(\alpha)} \, t^{\alpha-1} &:1<\alpha<2,\\
   			\frac{\bar G_{\infty}}{c}\, \frac{t}{\log t} &:\alpha=2,\\
   			\frac{\bar G_{\infty}}{\bar H_{\infty}}\, t &:\alpha>2.
          		\end{cases}
		\end{eqnarray*}
	\end{itemize}
	\item $\varrho=0$
			\begin{eqnarray*}
			\E[u(t,k)^2] \sim \begin{cases}
			\frac{\kappa c}{1-\alpha}t^{1-\alpha}&:\alpha<1,\\
   			\kappa c \log(t)&:\alpha=1,\\
   			1+\kappa \bar G_{\infty}&:\alpha>1,
          		\end{cases}\quad\text{ as }t\to \infty.
          		\end{eqnarray*}
	\item $\varrho<0$
	\begin{eqnarray*}
		\E[u(t,k)^2] \sim 
		\begin{cases}
			1-\frac{1}{\varrho}  &:\alpha\leq 1,\\
			1-\frac{1}{\varrho}+\frac{1}{\varrho(1-\varrho\kappa \bar G_{\infty})}&:\alpha>1,
		\end{cases}\quad\text{ as }t\to \infty.
	\end{eqnarray*}
	\end{itemize}
\end{prop}
For $\varrho=1$ (parabolic Anderson model), the subexponential growth was partially analysed for finite range transitions in \cite{DD07} (see their page 15). An example which was not included is for instance the Riemann walk defined in Example \ref{e1}. Since in this case
\begin{align}
	\bar p_t(0,0)\sim c t^{-1/\beta},\quad\text{as }t\rightarrow \infty,
\end{align}
it serves as a convenient example for the above results which exhibits a precise recurrence/transience transition at $\beta=1$. Further, the simple random walk on $\Z^d$ is contained with
\begin{eqnarray}\label{rates:laplace}
	\bar p_t(0,0)\sim c t^{-d/2},\quad\text{as }t\rightarrow \infty.
\end{eqnarray}

Combining the intermittency result with the extension of the results of \cite{CK00} given in \cite{BDE09}, we support the unstable behaviour of symbiotic branching for $\varrho>0$. It is quite standard (see for instance \cite{GM90}) that spatial processes being intermittent have a very local property: For large times the mass of the process is concentrated on few sites (``islands'').  Since for symbiotic branching, on each finite box solutions approach each constant configuration infinitely often, the islands do not stabilize. Since the diffusion function has the form $\sqrt{\kappa u(t,k)v(t,k)}$, we see that $u$ will not produce high peaks if $v$ is very small and vice versa. Hence, we suspect that $u,v$ are concentrated on the same islands. Understanding the pathwise behaviour better is an ambitious task for the future.
\subsection{Aging}
Recently in \cite{DD07} the concept of aging was discussed for certain classes of interacting diffusions. They say that aging (for linear test-functions) appears if the limit
\begin{eqnarray*}
	\lim_{t,s\rightarrow \infty}\corr[u(t,k),u(t+s,k)]
\end{eqnarray*}
depends on the choices of $s$ and $t$. Aging does not appear if this is not the case. The main results of \cite{DD07} were formulated with more general test-functions, though, restricted to finite range transitions. Differently, the present technique is restricted to linear test-functions but not to finite range transition. Our results suggest that neither finite range nor the linearity of test-functions is crucial. Symmetry of the transitions is assumed as in \cite{DD07}.\\

In \cite{DD07} it is shown that no aging appears in the parabolic Anderson model (in our model $\varrho=1$) in any dimensions for the discrete Laplacian. Further, for the super random walk (in our model related to $\varrho=0$) it was shown that aging appears exactly in dimensions $1,2$. This leads to the question if there are different phases for the symbiotic branching model. We show that the model exhibits three different regimes; an Anderson model like behaviour for $\varrho>0$, a super random walk like behaviour for $\varrho=0$, and a stepping stone model like behaviour for $\varrho<0$. The new case $\varrho<0$ and Corollary \ref{coro} suggest that there are three regimes in which the most prominent examples fall.
\begin{thm}[Aging for Symbiotic Branching]\label{t5} Let $(u_t,v_t)$ be a solution of the symbiotic branching model with homogeneous initial conditions and $a(i,j)=a(j,i)$. Then, if $\bar p_t(0,0)\sim c  t^{-\alpha}$, as $t\rightarrow \infty$, the following is true.
	\begin{itemize}
		\item[i)] 
			If $\varrho>0$, then no aging occurs for any $\alpha>0$.
		\item[ii)] If $\varrho=0$, then
			\begin{itemize}		
				\item no aging occurs, for any $\alpha>1$,
				\item $\lim_{t,s\rightarrow \infty,\log(s)/\log(t)=a}\corr[u(t,k),u(t+s,k)]=(1-a)_+$, for $\alpha=1$,
				\item $\lim_{t,s\rightarrow \infty,s=at}\corr[u(t,k),u(t+s,k)]=\frac{(1+\frac{a}{2})^{1-\alpha}-(\frac{a}{2})^{1-\alpha}}{(1+a)^{\frac{1-\alpha}{2}}}$, for any $\alpha<1$.
			\end{itemize}
		\item[iii)] If $\varrho<0$, then
			\begin{itemize}		
				\item no aging occurs, for any $\alpha>1$,
				\item $\lim_{t,s\rightarrow \infty,\log(s)/\log(t)=a}\corr[u(t,k),u(t+s,k)]=(1-a)_+$, for $\alpha=1$,
				\item $\lim_{t,s\rightarrow \infty,s=at}\corr[u(t,k),u(t+s,k)]= \frac{\int_0^1 (2 r+ a)^{-\alpha} (1-r)^{\alpha-1} d r} {
2^{-\alpha} \Gamma(\alpha)\Gamma(1-\alpha)}$, for any $\alpha<1$.
			\end{itemize}
	\end{itemize}
\end{thm}

We emphasise that our proof of Theorem \ref{t5} can be applied to more general interacting diffusions. In particular, the examples from the introduction are included. For finite range transitions, i), ii), iii) of the following proposition were proven in \cite{DD07}. 
\begin{prop}\label{coro}Consider solutions of (\ref{111}) with homogeneous initial conditions.  Then for
	\begin{itemize}
		\item[i)] $f(u)=u^2$, aging appears as in Theorem \ref{t5} i),
		\item[ii)] $0<\alpha_1\leq f(u)\leq \alpha_2$, aging appears as  in Theorem \ref{t5} ii),
		\item[iii)] $f(u)=u$, aging appears as in Theorem \ref{t5} ii),
		\item[iv)] $f(u)=u(1-u)$, aging appears as in Theorem \ref{t5} iii).
	\end{itemize}
	In the cases in which aging occurs, the upper and lower limits are bounded by the stated values up to constants depending on $f$.
\end{prop}

\section{Proofs}\label{section2}
\subsection{Some Results on Symbiotic Branching}\label{sm}
We start with a discussion on how second moments of symbiotic branching processes can be reduced to exponential moments of local times. Let us first recall the two-types particle moment-dual introduced in Section 3.1 of \cite{EF04}. This will be used to calculate second moments explicitly. Since the dual Markov process is described formally in \cite{EF04} we only sketch the pathwise behaviour. To find a suitable description of the mixed moment $\E[u(t,k)^nv(t,k)^m]$, $n+m$ particles are located at position $k\in \Z^d$. Each particle moves independently as a continuous time Markov process on $\Z^d$ with transition rates given by $(a(i,j))_{i,j\in\Z^d}$. At time $0$, $n$ particles have type $1$, $m$ particles have type $2$. One particle of each pair changes its type when the time the two particles have spent at same sites with same type exceeds an independent exponential time with parameter $\kappa$. Let
\begin{eqnarray*}
	L_t^=&=&\text{total collision time of all pairs of same type up to time }t,\\
	L_t^{\neq}&=&\text{total collision time of all pairs of different type up to time }t,\\
	l^1_t(a)&=&\text{number of particles of type }1\text{ at site }a\text{ at time t},\\
	l^2_t(a)&=&\text{number of particles of type }2\text{ at site }a\text{ at time t},\\
	(u_0,v_0)^{l_t}&=&\prod_{a\in \Z^d}u_0(a)^{l_t^1(a)}v_0(a)^{l_t^2(a)}.
\end{eqnarray*}
Note that since there are only $n+m$ many particles the infinite product is actually a finite product and hence well-defined. 
\begin{lem}\label{mom}
	Let $(u_t,v_t)$ be a solution of the symbiotic branching model with initial conditions $u_0,v_0\in E$ and $\varrho\in[-1,1]$. Then, for any $k\in\Z^d$, $t\geq 0$,
	\begin{eqnarray*}
		\E[u(t,k)^nv(t,k)^m]&=&\E\big[(u_0,v_0)^{l_t}e^{\kappa(L_t^=+\varrho L_t^{\neq})}\big].
	\end{eqnarray*}
\end{lem}
Though the proof for the moment-duality was given in \cite{EF04} (see the proof of their Proposition 9) only for the discrete Laplacian we skip a proof. For general transitions $(a(i,j))_{i,j\in\Z^d}$ the proof follows along the same lines.\\
Note that for homogeneous initial conditions $u_0=v_0\equiv 1$ the first factor in the expectation of the right-hand side equals $1$. Lemma \ref{mom} in the special case $\varrho=1$, $u_0=v_0\equiv 1$ was already stated in \cite{CM94}, reproven in \cite{GdH07}, and used to analyse the Lyapunov exponents of the parabolic Anderson model.\\
For $\varrho\neq 1$, the difficulty of the dual process is based on the two stochastic effects: On the one hand, one has to deal with collision times of random walks which were analysed in \cite{GdH07}. Additionally, particles have types either $1$ or $2$ which change dynamically.\\
Second moments are special since particles of different types do not change types anymore. Hence, when starting with two particles of same type there is precisely one event of changing types. This is used to obtain the following representation of second moments.
\begin{lem}\label{duality}
	Let $(u_t,v_t)$ be a solution of the symbiotic branching model with homogeneous initial conditions. Then, for any $k\in\Z^d$, $t\geq 0$,
	\begin{align*}
		\E[u(t,k)v(t,k)]&=\E[e^{\kappa\varrho L_t}],\\
		\E[u(t,k)^2]&=\E[v(t,k)^2]=\begin{cases}
		               1+\kappa\E[L_t]&:\varrho=0,\\
				1-\frac{1}{\varrho}+\frac{1}{\varrho}\E[e^{\kappa\varrho L_t}]&:\varrho\neq 0,
		               \end{cases}
		\end{align*}
	where $L_t$ denotes the local time in $0$ of the symmetrization $(\bar X_t)$ defined in (\ref{sy}) started in $0$.
\end{lem}
\begin{proof}
	The first expression for the mixed second moment follows directly from Lemma \ref{mom}: There are two particles which start with different types. Since pairs of particles of different types are never forced to change their types, they stay of different type for all time. Hence, $L_t^==0, L_t^{\neq}=L_t$ for all $t\geq 0$ and the assertion follows.\\
	For the second expression note that there is only one possible change of types. Starting with two particles of same types one of the types may change and the particles can not change their types again. Using independence of the particles and the exponential time we can make this explicit. Let $Y$ be an exponential variable with parameter $\kappa$, denote by $X$ the law of the two independent Markov processes, and $L_t$ their collision local time. Integrating out the exponential variable leads to
	\begin{align*}
		\E[u(t,k)^2]
		&=\E^{X\times Y}[e^{\kappa(L_t^=+\varrho L_t^{\neq})}]\\
		&=\E^{X\times Y}[e^{\kappa(L_t^=+\varrho L_t^{\neq})}1_{Y<L_t}]+\E^{X\times Y}[e^{\kappa(L_t^=+\varrho L_t^{\neq})}1_{Y\geq L_t}]\\
		&=\E^{X}\Big[\int_0^{L_t}\kappa e^{-\kappa x}e^{\kappa x+\kappa \varrho (L_t-x)}\,dx\Big]+\E^X\big[e^{\kappa L_t}\E^Y[1_{Y\geq L_t}]\big]\\
		&=\begin{cases}
		   \kappa\E[L_t]+\E[e^{\kappa L_t}e^{-\kappa L_t}]&:\varrho=0,\\
		   \E\Big[e^{\kappa\varrho L_t}\int_0^{L_t}\kappa e^{-\kappa \varrho x}\,dx\Big]+\E[e^{\kappa L_t}e^{-\kappa L_t}]&:\varrho\neq 0.
		  \end{cases}	
	\end{align*}		
	This proves the assertion.
\end{proof}
Now we prepare for the proof of the aging result.

\begin{lem}\label{appli}
	Let $(u_t,v_t)$ be a solution of the symbiotic branching model with homogeneous initial conditions and symmetric transitions $(a(i,j))_{i,j\in\Z^d}$. Then, for any $k\in\Z^d$, $t\geq 0$,
	\begin{eqnarray*}
		\E[u(t,k)u(t+s,k)]=1+\kappa\int_0^tp_{2r+s}(k,k)\E[e^{\kappa \varrho L_{t-r}}]\,dr
	\end{eqnarray*}
	and similarly for $v$.
\end{lem}
\begin{proof}
	The proof is only given for $u$ since due to symmetry the same proof works for $v$. We first employ the standard pointwise representation of solutions
	\begin{align}\label{gr}
		u(t,k)&=1+\sum_{i\in \Z^d}\int_0^tp_{t-s}(i,k)\sqrt{\kappa u(s,i)v(s,i)}\,dW^1_s(i)
	\end{align}	
	yielding
	\begin{align*}
		\E[&u(t,k)u(t+s,k)]\\
		&=1+\E\Big[\sum_{i\in\Z^d}\int_0^tp_{t-r}(i,k)\sqrt{\kappa u(r,i)v(r,i)}\,dW^1_r(i)
		\sum_{j\in\Z^d}\int_0^{t+s}p_{t+s-l}(j,k)\sqrt{\kappa u(l,j)v(l,j)}\,dW^1_l(j)\Big].
	\end{align*}
	Further, since martingale increments are orthogonal this equals
	\begin{align*}		
		1+\E\Big[\sum_{i\in\Z^d}\int_0^tp_{t-r}(i,k)\sqrt{\kappa u(r,i)v(r,i)}\,dW^1_r(i)
		\sum_{j\in\Z^d}\int_0^{t}p_{t+s-l}(j,k)\sqrt{\kappa u(l,j)v(l,j)}\,dW^1_l(j)\Big].
	\end{align*}	
	Now using independence of $W^1(i),W^1(j)$ for $i\neq j$ and It\^o's isometry we continue the chain of equalities as
	\begin{align*}
		1&+\sum_{i\in\Z^d}\E\left[\int_0^tp_{t-r}(i,k)p_{t+s-r}(i,k)\kappa u(r,i)v(r,i)\,dr\right]\\
		&=1+\int_0^t\sum_{i\in\Z^d} p_{t-r}(i,k)p_{t+s-r}(i,k)\kappa\E[u(r,i)v(r,i)]\,dr,
	\end{align*}	
	where we were allowed to change the order of integration since all terms are non-negative. Using Lemma \ref{duality}, which in particular shows for homogeneous initial conditions that second moments do not depend on the spatial variable, symmetry of the transitions, and the Chapman-Kolmogorov equality, we finish with
	\begin{align*}
		1&+\int_0^t\sum_{i\in\Z^d} p_{t-r}(k,i)p_{t+s-r}(i,k)\kappa \E[e^{\kappa \varrho L_{r}}]\,dr\\
		&=1+\kappa\int_0^tp_{2t+s-2r}(k,k) \E[e^{\kappa \varrho L_{r}}]\,dr
		=1+\kappa\int_0^tp_{2r+s}(k,k) \E[e^{\kappa \varrho L_{t-r}}]\,dr.
	\end{align*}
\end{proof}
Since we are going to examine the second Lyapunov exponent $\gamma_2$ of solutions we give a simple argument which ensures existence of the exponent.
\begin{lem}\label{existencegamma}
	Let $(u_t,v_t)$ be a solution of the symbiotic branching model with homogeneous initial conditions. Then the Lyapunov exponent $\gamma_2$ exists.
\end{lem}
\begin{proof}
	Note that to ensure existence of the limits 
	\begin{eqnarray*}
		\lim_{t\rightarrow \infty}\frac{1}{t}\log\E[u(t,k)^2]
	\end{eqnarray*}
	it suffices to show subadditivity of $\log \E[u(t,k)^2]$. Using Lemma \ref{duality} this is reduced to showing subadditivity of $\log\E[e^{\kappa \varrho L_t}]$, where $L_t$ is the local time in $0$ of $(\bar X_t)$ started in $0$. Thus, by conditioning on $\bar X_s$, we get
	\begin{align*}
		\log \E^0[e^{\kappa\varrho L_{t+s}}]
		=\log \E^0[e^{\kappa\varrho L_{s}}\E^{\bar X_s}[e^{\kappa\varrho L_t}]]
		\leq\log \E^0[e^{\kappa\varrho L_{s}}\E^{0}[e^{\kappa \varrho L_t}]]
		=\log \E^0[e^{\kappa\varrho L_{s}}]+\log\E^{0}[e^{\kappa \varrho L_t}].
	\end{align*}
\end{proof}

\subsection{Exponential Moments of Local Times}\label{sectionlocal}
\subsubsection{Preliminaries}
In Lemma~\ref{duality}, we observed that in order to study second moments of symbiotic branching processes it suffices to study exponential moments of local times of the symmetrization $(\bar X_t)$. We now take up this issue and discuss exponential moments of $L_t$ in greater generality than needed for the symbiotic branching model. For the following let $(X_t)$ be a time-homogeneous Markov process with countable state space $S$ and transition kernel $(a(i,j))_{i,j\in S}$. In particular, the transition rates are not assumed to be symmetric. \\
We start with a renewal-type equation for exponential moments of local times.
\begin{lem}\label{l3}Let $L_t$ be the local time of $(X_t)$ in $i\in S$ for the process started in $i$. Then for $\kappa\in\R$ the following equation holds:
	\begin{eqnarray}
		\E[e^{\kappa L_t}]=1+\kappa \int_0^tp_r(i,i)\E[e^{\kappa L_{t-r}}]\,dr,\quad t\geq 0 \label{eqn:maineqn}.
	\end{eqnarray}	
\end{lem}
\begin{proof}
	We use the exponential series to get
	\begin{align*}
		\E[e^{\kappa L_t}]
                &=\E\Big[e^{\kappa \int_0^t \delta_i(X_s)\,ds}\Big]
		=\E\left[\sum_{n=0}^{\infty}\frac{\kappa^n}{n!}\left(\int_0^t \delta_i(X_{s})\,ds\right)^n\right]\\
		&=1+\E\left[\sum_{n=1}^{\infty}\frac{\kappa^n}{n!}\int_0^t \cdots \int_0^t \delta_i(X_{s_1})\ldots \delta_i(X_{s_n})\,ds_n\ldots ds_1\right]\\
		&=1+\E\left[\sum_{n=1}^{\infty}\kappa^n\int_0^t \int_{s_1}^{t} \cdots \int_{s_{n-1}}^{t} \delta_i(X_{s_1})\ldots \delta_i(X_{s_n})\,ds_n\ldots d s_2 ds_1\right].
	\end{align*}
	The last step is justified by the fact that the function that is integrated is symmetric in all arguments and, thus, it suffices to integrate over a simplex. We can exchange sum and expectation and obtain that the last expression equals
	\begin{eqnarray*}
		&&1+ \kappa \int_0^t \sum_{n=1}^{\infty}\kappa^{n-1} \int_{s_1}^{t} \cdots \int_{s_{n-1}}^{t} \P [ X_{s_1}=i, \ldots, X_{s_n}=i]\,ds_n\ldots d s_2 ds_1.
	\end{eqnarray*}
	Due to the Markov property, the last expression equals
	\begin{align*}
		1&+ \kappa \int_0^t p_{s_1}(i,i) \sum_{n=1}^{\infty}\kappa^{n-1} \int_{s_1}^{t} \cdots \int_{s_{n-1}}^{t} \P [ X_{s_2-s_1}=0, \ldots, X_{s_n-s_1}=0] \,ds_n\ldots d s_2 ds_1
	\end{align*}
	and can be rewritten as
	\begin{align*}
		1+ \kappa \int_0^t p_{s_1}(i,i) \left(\sum_{n=1}^{\infty}\kappa^{n-1} \int_{0}^{t-s_1} \cdots \int_{s_{n-1}}^{t-s_1} \P [ X_{s_2}=0, \ldots, X_{s_n}=0] \,ds_n\ldots d s_2 \right) ds_1.
	\end{align*}
	Using the same line of arguments backwards for the term in parenthesis, the assertion follows.
\end{proof}

\begin{rem}
	A similar renewal-type equation as (\ref{eqn:maineqn}) can be shown with essentially the same proof for a discrete-time Markov process. It reads
	\begin{eqnarray*}
		\E[e^{\kappa L_m}] = 1 + \kappa \sum_{n=0}^{m} p_n(i,i) \E[e^{\kappa L_{m-n}}],\quad m\geq 1,
	\end{eqnarray*}
	where $p_n(i,i)$ is the return probability after $n$ steps and $L_n$ is the number of visits after $n$ steps.\\
	Similar equations were obtained for symmetric Markov chains on $\Z^d$ in \cite{MR92} using a completely different technique. Note that neither symmetry nor any structure of the set $S$ is needed. The information on the geometry of $S$ is completely encoded in $p_t(i,i)$. 
\end{rem}
	For the rest of this section we fix the Markov process $(X_t)$, $i\in S$, and abbreviate
	\begin{eqnarray*}
	 f(t)=p_t(i,i),\quad g(t)=\E[e^{\kappa L_t}].
	\end{eqnarray*}
	The return probabilities $p_t(i,i)$ are always assumed to be strongly asymptotically equivalent to $ct^{-\alpha}$, as $t \to\infty$, for $\alpha>0$, as for instance for simple random walks on $\Z^d$ and the Riemann walk on $\Z$. Further, $f$ is monotone, decreasing, positive with $f(0)=1$, and $g$ is monotone, increasing, positive with $g(0)=1$. The Laplace transform for a function $h$ on $\R_{\geq 0}$ is denoted by $\hat h$ and the convolution of two functions $f,g$ is denoted $f\ast g$. In this notation Equation \ref{eqn:maineqn} reads
	\begin{align}\label{db}
		g(t)=1+\kappa (f\ast g)(t), \quad t\geq 0.
	\end{align}
	Taking the Laplace transform of Equation (\ref{db}) leads to
	\begin{align} \label{eqn:kappa0prime}
 		\hat{g}(\lambda) = \frac{1}{\lambda} + \kappa \hat{f}(\lambda) \hat{g} (\lambda),\quad \lambda>0.
	\end{align}
	Obviously, since $f$ is bounded by $1$, $\hat{f}(\lambda)$ is always finite for all $\lambda\geq 0$. A priori this is not true for $g$ but if so, we obtain a useful representation from (\ref{eqn:kappa0prime}).
			\begin{lemma}\label{le}
				If $\hat{g}(\lambda)<\infty$, then
				\begin{eqnarray}
					\hat{g}(\lambda) = \frac{1}{\lambda ( 1 - \kappa \hat{f}(\lambda))}.
				\end{eqnarray}
			\end{lemma}
			In the following we proceed in two steps. First, we use (\ref{db}) to understand in which cases $g(t)$ grows exponentially in $t$ and discuss properties of the exponential growth rate. The following corresponence between exponential growth and finiteness of Laplace transforms holds (existence of the limit was proven in Lemma \ref{existencegamma}):
			\begin{eqnarray}\label{a}
				 \lim_{t\rightarrow \infty}\frac{1}{t}\log g(t)\geq c\quad \textrm{if and only if} \quad \hat{g}(c)=\infty.	
			\end{eqnarray}
			This observation is particularly important for the second step in which we discuss the behaviour of $g(t)$ as $t\to\infty$. In the cases in which $g(t)$ grows subexponentially (\ref{a}) implies that $\hat{g}(\lambda)<\infty$ for all $\lambda>0$. Hence, Lemma \ref{le} can be used for all $\lambda>0$. The strategy in this case is the following: By assumption, the asymptotic behaviour of $f(t)$ as $t$ tends to infinity is known, namely $ct^{-\alpha}$. Using Tauberian theorems the asymptotic behaviour of $\hat{f}(\lambda)$ as $\lambda$ tends to zero can be deduced. By Lemma \ref{le} this determines the asymptotic 	behaviour of $\hat{g}(\lambda)$ as $\lambda$ tends to zero. Using Tauberian theorems in the opposite direction, the asymptotic behaviour of $g(t)$ as $t$ tends to infinity is obtained.\\
			
			To manage the transfer from the behaviour of $f$ to $\hat{f}$ and back from $\hat{g}$  to $g$ the following Tauberian theorems are used. They are taken from \cite{BGT89} (see Theorem~1.7.6, Theorem 1.7.1, Corollary~8.1.7, and the considerations at the beginning of Section 8.1, \S 3).
			\begin{lemma} \label{lem:moment-duality} Let $h$ be a monotone function on $\R_{\geq0}$ with $h(0)=1$, then the following hold:
				\begin{enumerate}
 					\item [i)] If $\alpha<1$ and $\delta\in\R$, then $h(t)\sim c t^{-\alpha} (\log t)^\delta\text{ as }t\rightarrow \infty$ if and only if
					\begin{eqnarray*}
						\hat{h}(\lambda)\sim c \Gamma(1-\alpha) \lambda^{\alpha-1} (\log(1/\lambda))^{\delta},
					\end{eqnarray*}
					as $\lambda\to 0$.
					\item[ii)] If $h(t)\sim c t^{-1}\text{ as }t\rightarrow \infty$, then 
					\begin{eqnarray*}
						\hat{h}(\lambda)\sim c  \log(1/\lambda),
					\end{eqnarray*}
					as $\lambda\to 0$.
					\item[iii)]
					If $\alpha>1$ and $h(t)\sim c t^{-\alpha}\text{ as }t\rightarrow \infty$, then $I:=\int_0^\infty h(t)\, d t<\infty$ and
					\begin{eqnarray*}
						I - \hat h(\lambda) \sim
						\begin{cases}
   							\frac{c \Gamma(2-\alpha)}{\alpha-1} \lambda^{\alpha-1} &: 1<\alpha<2,\\
   							c \lambda \log \big(\frac{1}{\lambda}\big) &: \alpha=2,\\
   						\lambda \int_0^\infty t\, h(t) \, dt &: \alpha>2,
						\end{cases}
					\end{eqnarray*}
					as $\lambda\to 0$.
				\end{enumerate} 
			\end{lemma}

\subsubsection{Analysis of $\kappa>0$, Exponential Growth}
			The main point of the analysis is the following representation of the exponential growth rate which follows directly from Lemma \ref{l3}.
			\begin{prop}\label{prop1}
				Let $\kappa>0$, then
				\begin{equation} \label{eqn:defr}
					r(\kappa):=\lim_{t\to \infty} \frac{1}{t}\log \E[e^{\kappa L_t}]=\hat{f}^{-1}\Big(\frac{1}{\kappa}\Big).
				\end{equation}
			\end{prop}
			\begin{proof} 
				First, (\ref{a}) implies that
				\begin{eqnarray*}
					\inf\{\lambda|\hat{g}(\lambda)<\infty\}=\lim_{t\rightarrow\infty}\frac{1}{t}\log g(t).
				\end{eqnarray*}
				Moreover, 
				\begin{eqnarray*}
					\hat{f}^{-1}\Big(\frac{1}{\kappa}\Big)=\inf\Big\{\lambda\Big|\hat{f}(\lambda)< \frac{1}{\kappa}\Big\}.
				\end{eqnarray*}
				We are done if we can show
				\begin{eqnarray*}
				\{\lambda|\hat{g}(\lambda)<\infty\}=\Big\{\lambda\Big|\hat{f}(\lambda)< \frac{1}{\kappa}\Big\}.
				\end{eqnarray*}
				First we show ``$\subseteq$''. Due to Lemma \ref{l3} we obtain $\hat{g}(\lambda)=\frac{1}{\lambda}+\kappa\hat{f}(\lambda)\hat{g}(\lambda)$ which implies $\hat{g}(\lambda)> \kappa \hat{g}(\lambda)\hat{f}(\lambda)$. Since $\hat g(\lambda)<\infty$ this shows that $\hat{f}(\lambda)< \frac{1}{\kappa}$.\\
				 Now we show ``$\supseteq$''. First, iterating (\ref{db}) yields for fixed $n$
				 \begin{eqnarray*}
				 	g(t)=\sum_{i=0}^n\kappa^i(f^{\ast i}\ast 1)(t)+\kappa^{n+1}(f^{\ast(n+1)}\ast g)(t).
				 \end{eqnarray*}
				Using $f(t)\leq 1$ and $g(t)=\E[e^{\kappa L_t}]\leq e^{\kappa t}$ yields
				\begin{eqnarray*}
					\kappa^{n+1}(f^{\ast (n+1)}\ast g)(t)&=&\kappa^{n+1}\int_0^tf^{\ast (n+1)}(s)g(t-s)\,ds\\
						&\leq&\kappa^{n+1}\int_0^t\frac{s^{n}}{n!}e^{\kappa(t-s)}\,ds
						\leq\kappa\frac{(\kappa t)^{n}}{n!}\int_0^te^{\kappa(t-s)}\,ds\rightarrow 0,
				\end{eqnarray*}
				as $n\rightarrow \infty$. Hence, for fixed $t\geq 0$
				 \begin{eqnarray*}
				 	g(t)=\sum_{i=0}^{\infty}\kappa^i(f^{\ast i}\ast 1)(t).
				 \end{eqnarray*}
				Taking Laplace transforms we note that $\hat{g}(\lambda)$ is finite if and only if the Laplace transform of the right-hand side is finite. However, using Fubini's theorem (note that only $\kappa>0$ needs to be considered) we obtain
				\begin{eqnarray*}
					\widehat{\left(\sum_{i=0}^{\infty}\kappa^i(f^{\ast i}\ast 1)\right)}(\lambda)=\sum_{i=0}^{\infty}\kappa^i\widehat{f^{\ast i}\ast 1}(\lambda)
					=\frac{1}{\lambda}\sum_{i=0}^{\infty}(\kappa\hat{f}(\lambda))^i,
				\end{eqnarray*}
				which is finite since we assumed $\kappa \hat{f}(\lambda)<1$.
			\end{proof}
			In particular, the previous result shows that understanding $\hat{f}^{-1}$ suffices to understand the exponential growth rates of $\E[e^{\kappa L_t}]$. This is not difficult due to the following observation: $\hat{f}$ is a strictly decreasing, convex function with $\hat{f}(0)=G_{\infty}$. Hence, $\hat{f}^{-1}$ is a strictly decreasing, convex function with $\lim_{\lambda\rightarrow 0}\hat{f}^{-1}(\lambda)=\infty$ and $\hat{f}^{-1}(\lambda)=0$ if and only if $\lambda \geq G_{\infty}$. This implies that $\hat{f}^{-1}\big(\frac{1}{\lambda}\big)=0$ precisely for $\lambda\leq \frac{1}{G_{\infty}}$. This and more properties of the exponential growth rate are collected in the following corollary.
			\begin{cor}\label{cor1}
				Let $\kappa>0$ and $r(\kappa)=\lim_{t\rightarrow\infty}\frac{1}{t}\log \E[e^{\kappa L_t}]$. Then with $\kappa_{cr}:=\frac{1}{G_{\infty}}$ the following hold:
				\begin{itemize}
					\item[i)] $r(\kappa)\geq 0$ and $r(\kappa)>0$ if and only if $\kappa>\kappa_{cr}$,
 					\item[ii)] the function $\kappa \mapsto r(\kappa)$ is strictly convex for $\kappa>\kappa_{cr}$,
 					\item[iii)] $r(\kappa)\leq \kappa$ for all $\kappa$, and $\frac{r(\kappa)}{\kappa} \to 1$, as $\kappa\to\infty$,
 					\item[iv)] if $\alpha\leq 1$, then $\kappa_{cr}=0$ and, as $\kappa\to 0$,
					\begin{align*}
						r(\kappa) \sim \begin{cases}
   						\kappa^{\frac{1}{1-\alpha}}(c \Gamma(1-\alpha))^{\frac{1}{1-\alpha}} &: 0<\alpha<1,\\
   						\exp( -(c \kappa)^{-1} + o(\kappa^{-1})) &: \alpha=1,\\
						\end{cases}
					\end{align*}
					\item[v)] if $\alpha>1$, then $\kappa_{cr}>0$ and, as $\kappa\searrow \kappa_c$,
					\begin{align*}
						r(\kappa) \sim \begin{cases}
   						(\kappa-\kappa_c)^{\frac{1}{\alpha -1}}\left(\frac{ G_{\infty}^2 (\alpha-1) }{c \Gamma(2-\alpha)} \right)^{\frac{1}{\alpha-1}} &: 1<\alpha<2,\\
   						\frac{\kappa-\kappa_c}{\log \big(\frac{1}{\kappa-\kappa_c}\big)}\frac{G_{\infty}^2}{c}  &: \alpha=2,\\
   						 (\kappa-\kappa_c)\frac{G_{\infty}^2}{H_{\infty}} &: \alpha>2.
						\end{cases}
					\end{align*}
				\end{itemize}
			\end{cor}
			\begin{proof}
				Parts i) and ii) are proven as argued above the corollary.\\
				Since $f\leq 1$, the first part of iii) follows from
				 \begin{eqnarray*}
				      	\frac{1}{\kappa}=\hat{f}(r(\kappa))=\int_0^{\infty}e^{-r(\kappa)x}f(x)\,dx\leq \int_0^{\infty}e^{-r(\kappa) x}\,dx=\frac{1}{r(\kappa)}.
				     \end{eqnarray*}
				    Continuity of $f$ and $f(0)=1$ imply that  for $\epsilon>0$ there is $x_0(\epsilon)$ such that $f(x)\geq1-\epsilon$ for $x\leq x_0(\epsilon)$. Hence,
				     \begin{align*}
				     	\frac{1}{\kappa}&=\hat{f}(r(\kappa))=\int_0^{\infty}e^{-r(\kappa)x}f(x)\,dx\\
				     	&\geq (1-\epsilon)\int_0^{x_0(\epsilon)}e^{-r(\kappa)x}\,dx
				     	=(1-\epsilon)\frac{1}{r(\kappa)}\big(1-e^{-r(\kappa)x_0(\epsilon)}\big).
				     \end{align*}
					Since $r(\kappa)\rightarrow \infty$ for $\kappa\rightarrow \infty$ we obtain
					\begin{eqnarray*}
					 \liminf_{k\rightarrow \infty}\frac{r(\kappa)}{\kappa}\geq 1.
					\end{eqnarray*}
					The second part of iii) now follows since as well $\frac{r(\kappa)}{\kappa}\leq 1$ for all $\kappa>0$.\\
				Finally, for iv) and v) note that the asymptotics of $\hat{f}$ for $\lambda\rightarrow 0$ are known from Lemma \ref{lem:moment-duality}. This translates to $\hat{f}^{-1}(\lambda)$ and hence to $r(\kappa)=\hat{f}^{-1}\big(\frac{1}{\kappa}\big)$.
			\end{proof}
			\subsubsection{Analysis of $\kappa>0$, Subexponential Growth}
			
			So far, we have understood the behaviour of $\E[e^{\kappa L_t}]=g(t)$ as $t\to\infty$ for $\kappa>\frac{1}{G_{\infty}}$. In this case $g(t)$ grows exponentially and the behaviour of the exponential rates in $\kappa$ could be analysed. We now come to the case $\kappa\leq \frac{1}{G_{\infty}}$. First, if $G_{\infty}=\infty$, there is nothing to be done since the only appearing case is $\kappa=0$ which yields $g(t)=1$ for all $t\geq 0$. Hence, we can stick to $G_{\infty}<\infty$.
			
			\begin{prop}\label{prop10}
				Let $\kappa>0$ and $\kappa<\frac{1}{G_{\infty}}$. Then, as $t\to \infty$,
				\begin{eqnarray*}
					\E[e^{\kappa L_t}] \sim \frac{1}{1-\kappa G_{\infty}}.
				\end{eqnarray*}
			\end{prop}

			\begin{proof}
				Since $G_{\infty}<\infty$ we can apply part iii) of Lemma~\ref{lem:moment-duality}. Hence, $\hat{f}(\lambda)\to G_{\infty}$,  as $\lambda\to 0$. As discussed above, since $g(t)$ does not grow exponentially, $\hat{g}(\lambda)>0$ for all $\lambda>0$, and we can use Lemma \ref{le}. This implies 
				\begin{eqnarray*}
					\hat{g}(\lambda)\sim \lambda^{-1}\frac{1}{(1-\kappa G_{\infty})} ,
				\end{eqnarray*}
				as $\lambda\to 0$. Going backwards with  Lemma~\ref{lem:moment-duality}, part i), $\alpha=\delta=0$, the asymptotic of $g$ follows.
			\end{proof}

			\begin{prop}\label{prop11}
				Let $\kappa>0$ and $\kappa=\frac{1}{G_{\infty}}$. Then, as $t\to\infty$,
				\begin{eqnarray*}
					\E[e^{\kappa L_t}] \sim \begin{cases}
   					t^{\alpha-1}\frac{\alpha-1}{\kappa c \Gamma(2-\alpha)\Gamma(\alpha)}  &:1<\alpha<2,\\
   					\frac{t}{ \log t} \frac{1}{\kappa c}&:\alpha=2,\\
   					t\frac{1}{\kappa H_{\infty}} &:\alpha>2.
          				\end{cases}
				\end{eqnarray*}
			\end{prop}
			\begin{proof}
				Since $G_{\infty}<\infty$ we can apply Lemma \ref{lem:moment-duality}, part iii). Hence, $\hat{f}(\lambda)\sim G_{\infty}$, as $\lambda\to 0$. As discussed above, since $g(t)$ does not grow exponentially, $\hat{g}(\lambda)>0$ for all $\lambda>0$, and we can use Lemma \ref{le}. Since $\kappa G_{\infty}=1$, the denominator $(1-\kappa \hat{f}(\lambda))$ appearing in Lemma \ref{le} does not behave like a constant and we cannot apply part i) of Lemma \ref{lem:moment-duality} with $\alpha=\delta=0$. Instead we use Lemma \ref{lem:moment-duality}, part iii), to obtain
				\begin{align*}
					\hat{g}(\lambda) = \frac{1}{\lambda\kappa} \frac{1}{G_{\infty}- \hat{f}(\lambda)} \sim \frac{1}{\lambda\kappa}
					\begin{cases}
    					\frac{\alpha-1}{c \Gamma(2-\alpha)} \lambda^{1-\alpha} &:1<\alpha<2,\\
   					\frac{1}{c} \lambda^{-1} (\log 1/\lambda)^{-1} &:\alpha=2,\\
   					\lambda^{-1} \frac{1}{H_{\infty}} &:\alpha>2,
					\end{cases}
				\end{align*}
				as $\lambda\to 0$. This, by Lemma \ref{lem:moment-duality}, part i), implies the assertion.
			\end{proof}
%
		\subsubsection{Analysis of $\kappa<0$} \label{sec:kappa<0}
			We now investigate Equation (\ref{eqn:maineqn}) for $\kappa<0$.
			\begin{prop}\label{prop12}
		 		If $\kappa<0$, then, as $t\to \infty$,
				\begin{eqnarray*}
					\E[e^{\kappa L_t}] \sim \begin{cases}
           				\frac{1}{t^{1-\alpha}}\frac{1}{-\kappa c \Gamma(1-\alpha) \Gamma(\alpha)}  &:0<\alpha<1,\\
           				\frac{1}{\log t}\frac{1}{-\kappa c }&:\alpha=1,\\
           				\frac{1}{-\kappa G_{\infty}+1}&:\alpha>1.
          				\end{cases}
				\end{eqnarray*}
			\end{prop}

			\begin{proof}
				First note that for $\kappa<0$, $g(t)=\E[e^{\kappa L_t}]<1$ and hence for all $\lambda>0$, $\hat{g}(\lambda)<\infty$ which validates the use of Lemma \ref{le}. This implies
				\begin{eqnarray*}
				 	\hat{g}(\lambda)=  \frac{1}{\lambda( 1 - \kappa  \hat{f}(\lambda))} \sim \frac{1}{-\kappa} \lambda^{-1}\frac{1}{ \hat{f}(\lambda)},
				\end{eqnarray*}	
				as $\lambda\to 0$. Using Lemma~\ref{lem:moment-duality} in both directions returns the assertion.
			\end{proof}
	
	\subsection{Proofs of the Main Results}\label{sectionmain}
		\begin{proof}[Proofs of Theorem \ref{t3} and Proposition \ref{cor}]
			These follow from Lemma \ref{duality} and Corollary \ref{cor1}.
		\end{proof}

		\begin{proof}[Proof of Proposition \ref{pro}]
			This follows from Lemma \ref{duality} and Propositions \ref{prop10} and \ref{prop11}.
		\end{proof}
		
		\begin{proof}[Proof of Theorem \ref{t5}]
			For this proof we always denote by $\sim$ the strong asymptotics at infinity and we abbreviate $p_t=p_t(k,k)$. Lemma \ref{appli} implies
			\begin{eqnarray*}
				\textrm{cor}[u(t,k),u(t+s,k)]=\frac{\int_0^t p_{2r+s} \E[e^{\varrho \kappa L_{t-r}}] \, d r}{\sqrt{\int_0^t p_{2r} \E[e^{\varrho \kappa L_{t-r}}] \, d r \, \int_0^{t+s}
				p_{2r} \E[e^{\varrho \kappa L_{t+s-r}}] \, d r}}.
			\end{eqnarray*}
		\emph{Step 1, $\varrho=0$:} First, assume $\alpha>1$ which implies $\int_0^\infty p_{2 r} d r <\infty$. Since
			\begin{align*}
			 	\int_0^t p_{2 r+s} d r \approx  \int_0^t (2 r+s)^{-\alpha} d r \approx  \int_s^{t+s} r^{-\alpha} d r \leq \int_s^\infty r^{-\alpha} d r \stackrel{s\to\infty}{\to} 0,
			\end{align*}
			we obtain, independently of the choice of $t$ and $s$, 
			\begin{eqnarray*}
				\textrm{cor}[u(t,k),u(t+s,k)]=
				\frac{\int_0^t p_{2r+s} \, d r}{\sqrt{\int_0^t p_{2r} \, d r \, \int_0^{t+s} p_{2r}\, d r}}
				\stackrel{s,t\to\infty}{\to} 0.
			\end{eqnarray*}
			Here, we used $ f\approx g$ if $0<\liminf f/g\leq \limsup f/g<\infty$.
			We now come to the case $\alpha=1$, where we get
			\begin{align*}
				\int_0^t p_{2 r+s} d r \sim c \int_0^t (2 r+s+1)^{-1} d r = \frac{c}{2}\, \log \Big(\frac{2 t +s +1}{s+1}\Big).
			\end{align*}
			Therefore, we have
			\begin{align}
				\textrm{cor}[u(t,k),u(t+s,k)]\sim\frac{\log\big(\frac{2 t +s +1}{s+1}\big)}{\sqrt{\log (2t +1) \, \log( 2(t+s)+1)}}. \label{eqn:age1}
			\end{align}
			For $s=t^a$ with $a\leq 1$ this expression behaves asymptotically as
			\begin{align*}
				\frac{\log \big(t^{1-a}\big)}{\sqrt{\log (t) \, \log( t)}} = 1-a.
			\end{align*}
			On the other hand, for $s=t^a$ with $a\geq 1$  the term in (\ref{eqn:age1}) behaves asymptotically as
			\begin{align*}
				\frac{\log(1)}{\sqrt{\log (2(t+t^{a}) +1) \, \log( 2 (t+t^{a}) +1)}} \stackrel{s,t\rightarrow \infty}{\to} 0.
			\end{align*}
			Hence, for $\log(s)/\log(t) =a$, we obtain $\textrm{cor}[u(t,k),u(t+s,k)] \sim (1-a)_+.$\\
			Now suppose $\alpha<1$. Then
			\begin{align*}
			 	\int_0^t p_{2 r+s} d r \sim c \int_0^t (2 r+s+1)^{-\alpha} d r = \frac{c}{2}\, \frac{(2 t+s+1)^{1-\alpha} -(s+1)^{1-\alpha}}{1-\alpha}.
			\end{align*}
			Therefore, we have
			\begin{align*}
				\textrm{cor}[u(t,k),u(t+s,k)]\sim\frac{(2 t+s+1)^{1-\alpha} -(s+1)^{1-\alpha}}{\sqrt{((2 t+1)^{1-\alpha} -1) \, ((2 (t+s)+1)^{1-\alpha} -1)}}.
			\end{align*}
			For $s=at$ this behaves asymptotically as
			\begin{align*}
			\frac{(2 +a)^{1-\alpha} -a^{1-\alpha}}{\sqrt{2^{1-\alpha} \, (2 (1+a))^{1-\alpha}}} = \frac{(1 +a/2)^{1-\alpha}
			-(a/2)^{1-\alpha}}{(1+a)^{(1-\alpha)/2}}.
			\end{align*}
		\emph{Step 2, $\varrho<0$:} Let us first consider $\alpha>1$. Since $c_1 \leq \E e^{\varrho \kappa L_{t-r}} \leq c_2$ this case is exactly the same as $\varrho=0$, $\alpha>1$. 
			Now suppose $\alpha=1$. In this case we have by Proposition \ref{prop12}
			\begin{align*}
				\int_0^t p_{2 r+s} \E e^{\varrho \kappa L_{t-r}} d r &\sim \frac{c}{-\kappa \varrho c } \int_1^{t-e} \frac{1}{(2 r+1+s) \log (t-r)} d r 
				=
				\frac{c}{-\kappa \varrho c } \int_e^{t-1} \frac{1}{2 (t-r)+1+s} \frac{1}{ \log (r)} d r.
			\end{align*}
We use the scaling $s=t^a$ with $a<1$. Let $0<\theta<1$. The integral above can be split from $e$ to $\theta t$ and $\theta t$ to $t-1$. We treat the first integral and show that its order is less than $(\log (t))^{-1}$. First note that in the range of integration
\begin{align*}
\frac{1}{2 (t-e)+1+s} \leq \frac{1}{2 (t-r)+1+s}\leq \frac{1}{2 (t-\theta t)+1+s},
\end{align*}

Therefore,
\begin{align*}
\int_e^{\theta t}\frac{1}{2 (t-r)+1+s} \frac{1}{ \log (r)} d r \approx \frac{1}{t} \, \int_e^{\theta t}\frac{1}{ \log (r)} d r \approx
\frac{\theta}{\log (t)}.
\end{align*}

On the other hand, the second integral can be treated as follows. In its range of integration we have $\frac{1}{\log (t)} \leq \frac{1}{\log (r)}\leq\frac{1}{\log (\theta t)}.$ Therefore,
\begin{align*}
\int_{\theta t}^{t-1} \frac{1}{2 (t-r)+1+s} \frac{1}{ \log (r)} d r \sim \frac{1}{\log (t)} \int_{\theta t}^{t-1} \frac{1}{2 (t-r)+1+s}  d
r = \frac{1}{\log (t)}  \, \frac{1}{2} \log \left(\frac{2(t - \theta t)+s}{2+s}\right).
\end{align*}
Thus,
\begin{align*}
\textrm{cor}[u(t,k),u(t+s,k)]
 \sim \frac{\frac{1}{\log (t)}  \, \frac{1}{2} \log \left(\frac{2(1 - \theta) t+s}{2+s}\right)}{\sqrt{\frac{1}{\log (t)}  \, \frac{1}{2} \log ((1 -
\theta) t)\, \frac{1}{\log (t)}  \, \frac{1}{2} \log ((1 - \theta)(t+s))}} \sim \frac{ \log \left(\frac{ t}{t^a}\right)}{\sqrt{\log  (t)\log (t)}} = 1-a.
\end{align*}
Analogously, for case $a\geq 1$. Therefore, we get $\textrm{cor}[u(t,k),u(t+s,k)]\sim (1-a)_+,$ whenever $\log (s)/\log (t) = a$.\\
For $\varrho<0$ only $\alpha<1$ is left: Here, we have
\begin{align*}
\int_0^t p_{2 r+s} \E e^{\varrho \kappa L_{t-r}} d r \sim \frac{c}{-\kappa \varrho c \Gamma(\alpha)\Gamma(1-\alpha)} \int_0^t (2
r+1+s)^{-\alpha} (t-r)^{\alpha-1} d r
\end{align*}

We set $s=a t$. The integral can be rewritten as
\begin{align*}
\int_0^1 (2 r t+1+ at )^{-\alpha} (t-t r)^{\alpha-1} t d r \sim \int_0^1 (2 r+ a)^{-\alpha} (1-r)^{\alpha-1} d r.
\end{align*}

The same way one can see that
\begin{align*}
\int_0^t p_{2 r} \E[e^{\varrho \kappa L_{t-r}}] d r \sim \frac{c}{-\kappa \varrho c \Gamma(\alpha)\Gamma(1-\alpha)} \int_0^1 (2 r)^{-\alpha}
(1-r)^{\alpha-1} d r.
\end{align*}

Thus,
\begin{align*}
\textrm{cor}[u(t,k),u(t+s,k)]
  &\sim \frac{\int_0^1 (2 r+ a)^{-\alpha} (1-r)^{\alpha-1} d r}{\int_0^1 (2 r)^{-\alpha} (1-r)^{\alpha-1} d r}\\ &= \frac{\int_0^1 (2 r+
a)^{-\alpha} (1-r)^{\alpha-1} d r} { 2^{-\alpha} B(\alpha,1-\alpha)}= \frac{\int_0^1 (2 r+ a)^{-\alpha} (1-r)^{\alpha-1} d r} {
2^{-\alpha} \Gamma(\alpha)\Gamma(1-\alpha)},
\end{align*}
when  $s=a t \to \infty$. Here, $B$ denotes the Beta function.\\

	\emph{Step 3, $\varrho>0$:} \def\grw{\lambda} The transient case ($\alpha>1$) with $\kappa\varrho<\frac{1}{\bar G_{\infty}}$ has already appeared in the case $\varrho=0$ for $\alpha>1$. $\E[e^{\varrho\kappa L_t}]$ is bounded due to Proposition \ref{prop10}. Hence, no aging occurs.\\
	For $\kappa\varrho>\frac{1}{\bar G_{\infty}}$, and for $\alpha\leq 1$, we proved in Corollary \ref{cor1} i) that $\E[e^{ \varrho \kappa L_t}]$ grows exponentially. This implies that there is a $\grw>0$ depending on the growth rate, i.e.\ on $\varrho$ and $\kappa$, (for example $\grw=\frac{1}{2}\hat{\bar p}^{-1}\big(\frac{1}{\varrho \kappa}\big)$ does the job) such that for all $t\geq 0$ and $0\leq r\leq t$ we have 
	\begin{eqnarray*}
		\E[e^{ \varrho \kappa L_{t-r}}] \leq e^{-\grw r} \E[e^{ \varrho \kappa L_t}].
	\end{eqnarray*}

	Therefore,
	\begin{eqnarray*}
		\int_0^t p_{2 r+s} \E[e^{ \varrho \kappa L_{t-r}}] \, d r \leq p_{s} \E[e^{ \varrho \kappa L_{t}}]  \int_0^t e^{-\grw r} \, d r \leq \frac{1}{\grw}\, p_{s} \E[e^{ \varrho \kappa L_{t}}].
	\end{eqnarray*}

	On the other hand, by the assumption $p_t\sim c t^{-\alpha}$, the renewal-type equation of Lemma \ref{l3}, and since $\varrho>0$,
	\begin{eqnarray*}
		\int_0^t p_{2 r} \E[e^{ \varrho \kappa L_{t-r}}] \, d r \geq c \int_0^t p_{r} \E[e^{ \varrho \kappa L_{t-r}}] \, d r =\frac{c}{\varrho \kappa}  \, \left( \E[e^{ \varrho \kappa L_{t}}] -1\right) \geq c' \E[e^{ \varrho \kappa L_{t} }].
	\end{eqnarray*}

	Putting these pieces together we obtain that
	\begin{eqnarray*}
		\textrm{cor}[u(t,k),u(t+s,k)] \leq c''\, \frac{p_{s} \E[ e^{ \varrho \kappa L_{t} }]}{\sqrt{\E[ e^{ \varrho \kappa L_{t} }]\E[ e^{ \varrho \kappa L_{t+s} }]}} = c''p_{s}\, \sqrt{ \frac{ \E[ e^{ \varrho \kappa L_{t} }]}{\E[ e^{ \varrho \kappa L_{t+s} }]}}. 
	\end{eqnarray*}
	Note that, since $\varrho>0$, the term with the square root is bounded by $1$. Since clearly $p_s$ tends to zero, the whole expression must tend to zero independently of how $t,s\to \infty$. \\
	The only case left is $\alpha>1$ and $\kappa\varrho=\frac{1}{\bar G_{\infty}}$. First we consider $1<\alpha<2$.
\begin{align*}
\int_0^t p_{2 r} \E[e^{ \varrho \kappa L_{t-r} }] \, d r \approx \int_0^t (r+1)^{-\alpha} (t-r)^{\alpha-1} \, d r = \int_0^1
(r+1/t)^{-\alpha} (1-r)^{\alpha-1} \, d r.
\end{align*}

This expression tends to infinity for $t\to\infty$. The rate is
\begin{align*}
\approx  \int_0^{1/2} (r+1/t)^{-\alpha} (1-r)^{\alpha-1} \, d r \approx \int_0^{1/2} (r+1/t)^{-\alpha} \, d r \approx (1/t)^{1-\alpha}
= t^{\alpha-1}.
\end{align*}

On the other hand,
\begin{align*}
\int_0^t p_{2 r+s} \E[e^{ \varrho \kappa L_{t-r} }] \, d r  \approx \int_0^t (r+s)^{-\alpha} (t-r)^{\alpha-1} \, d r =\int_0^1
(r+s/t)^{-\alpha} (1-r)^{\alpha-1} \, d r.
\end{align*}

This expression is bounded or tends to infinity (depening on how $t,s\to \infty$). It is bounded by
\begin{align*}
c \int_0^{1/2} (r+s/t)^{-\alpha} (1-r)^{\alpha-1} \, d r \approx \int_0^{1/2} (r+s/t)^{-\alpha} \, d r \leq c_1 + c_2
(s/t)^{1-\alpha}.
\end{align*}

Putting these pieces together we obtain that
\begin{align*}
\frac{\int_0^t p_{2 r+s} \E[ e^{ \varrho \kappa L_{t-r} }] \, d r }{\sqrt{\int_0^t p_{2 r} \E[ e^{ \varrho \kappa L_{t-r} }] \, d r \int_0^{t+s}
p_{2 r} \E[ e^{ \varrho \kappa L_{t+s-r} }] \, d r}} \leq c' \,\frac{ c_1 + c_2 (s/t)^{1-\alpha}}{\sqrt{t^{\alpha-1} (t+s)^{\alpha-1}}} \to
0.
\end{align*}

The calculation is completely analogous in the cases $\alpha\geq 2$.
\end{proof}
		
\begin{proof}[Proof of Proposition \ref{coro}]
	The pointwise representation of (\ref{gr}) is not restricted to the symbiotic branching model but also holds for the interacting diffusions of this proposition. Hence, the same derivation as for the symbiotic branching model yields
	\begin{align*}
			\corr[u_t(k),u_{t+s}(k)]
			= \frac{\int_0^tp_{2r+s}(k,k)\E[f(u(t-r,k))]\,dr}{\sqrt{\int_0^tp_{2r}(k,k)\E[f(u(t-r,k))]\,dr\int_0^{t+s}p_{2r}(k,k)\E[f(u(t+s-r,k))]\,dr}}.
	\end{align*}
	Part i) is contained in Theorem \ref{t3} since the parabolic Anderson model appears as special case $\varrho=1$. For part ii) we can estimate the expectations from above and below to obtain the same result as in Theorem \ref{t5} ii) except constants. The same is true for part iii) since the pointwise representation (\ref{gr}) implies $\E[u(t,k)]=1$. Finally, we could interprete part iv) as a submodel of the symbiotic branching model. Instead, we give a direct proof using use the coalescing particles dual of \cite{S88}. The dual process consists of two independent particles started in $k$, performing transitions $(a(i,j))_{i,j\in\Z^d}$ in continuous time. After spending an exponential time $Y$ with parameter $\kappa$, independent of the particles, at same sites, the particles coalesce. We denote $X$ the law of the particles and suppose $u_0\equiv w\in(0,1)$. Then
	\begin{align*}
		\E[u(t,k)^2]&=\E^{Y\times X}[w^{\text{number of non-coalesced particles}}]\\
		&=\E^{Y\times X}[w1_{Y\leq L_t}]+\E^{Y\times X}[w^21_{Y> L_t}]
		=w(1-\E^X[e^{-\kappa L_t}])+w^2\E^X[e^{-\kappa L_t}],
	\end{align*}
	where $L_t$ denotes the collision time of the particles. Using $\E[u(t,k)]=w$ this yields $\E[f(u(t,k))]=\E[u(t,k)]-\E[u(t,k)^2]=(w-w^2)\E[e^{-\kappa L_t}]$ and we can proceed as for the symbiotic branching model with $\varrho<0$.
	\end{proof}

\section*{Acknowledgement}
	We are grateful to J.\ G\"artner from TU Berlin for simplifying our original proof of Lemma \ref{l3}.

\end{document}